\documentclass[12pt]{amsart}

\usepackage{amsmath, amscd, graphicx, latexsym, hyperref, rlepsf, times, }

\oddsidemargin 0in \theoremstyle{plain} \topmargin 0in

\oddsidemargin .25in \theoremstyle{plain}

\newtheorem{Thm}{Theorem}[section]

\newtheorem{Def}[Thm]{Definition}
\newtheorem{Rem}[Thm]{Remark}

\def\p{\partial}

\def\v{\vskip.12in}
\def\a{\alpha}

\def\e{\epsilon}
\def\d{\delta}
\def\G{\Gamma}
\def\g{\gamma}

\def\S{\Sigma}

\begin{document}

\title{Contact handle decompositions}


\author{Burak Ozbagci}

\begin{abstract}

We review Giroux's contact handles and contact handle attachments
in dimension three and show that a bypass attachment consists of a
pair of \emph{contact} 1 and 2-handles. As an application we
describe explicit contact handle decompositions of infinitely many
pairwise non-isotopic overtwisted $3$-spheres. We also give an
alternative proof of the fact that every compact contact
$3$-manifold (closed or with convex boundary) admits a contact
handle decomposition, which is a result originally due to Giroux.

\end{abstract}

\address{Department of Mathematics \\ Ko\c{c} University \\ Istanbul, Turkey}

\email{bozbagci@ku.edu.tr} 


\thanks{}

\v \v \v

\maketitle

\setcounter{section}{-1}


\section{Introduction}
Emmanuel Giroux announced the following result in a series of
lectures he delivered at Stanford University in the year 2000:
``Every contact $3$-manifold is convex" --- which signified the
closure of the program he initiated in his convexity paper
published in 1991, where he proved that every oriented
$3$-manifold has \emph{some} convex contact structure (\cite{gi},
Theorem III. 1.2). Apparently, an essential motivating factor for
studying convexity in contact topology is the following
straightforward consequence of the convexity theorem:
``\emph{Every contact $3$-manifold (closed or with convex
boundary) admits a contact handle decomposition}". We should point
out that for a closed contact $3$-manifold the existence of a
contact handle decomposition and the existence of an adapted open
book decomposition are equivalent. Despite the fact that several
explicit examples of adapted open book decompositions of closed
contact $3$-manifolds have been published and fruitfully used in
various other constructions since Giroux's breakthrough in 2000,
explicit examples of contact handle decompositions of closed
contact $3$-manifolds have not yet appeared in the literature. In
this article, we show that a bypass attachment  \cite{hon}
consists of a (topologically cancelling) pair of \emph{contact} 1
and 2-handles. As an application, for each positive integer $n$,
we describe an explicit contact handle decomposition of the
overtwisted $3$-sphere whose $d_3$-invariant is $(2n+1)/2 $.
Recall that
 two overtwisted contact structures are isotopic if and
only if they are homotopic as oriented $2$-plane fields
\cite{eli}. Moreover the homotopy classes of oriented $2$-plane
fields on $S^3$ are classified by their $d_3$-invariants (see
\cite{G} or \cite{ozst} for a detailed discussion).

For the sake of completeness, we also offer a alternative proof of
Giroux's handle decomposition theorem for compact contact
$3$-manifolds (closed or with convex boundary). Our proof is based
on a recent result due to Honda, Kazez and Mati\'{c} (\cite{hkm},
Theorem 1.1), which asserts that every compact contact
$3$-manifold with convex boundary has an adapted partial open book
decomposition. The technique that Honda, Kazez and Mati\'{c} apply
in constructing adapted partial open book decompositions of
contact $3$-manifolds with convex boundary is a generalization of
Giroux's method of constructing adapted open book decompositions
of closed contact $3$-manifolds. Giroux's construction, in turn,
is based on contact cell decompositions of contact $3$-manifolds
\cite{gir}. Hence the existence of contact \emph{handle}
decompositions of compact contact $3$-manifolds can be viewed as a
consequence of the existence of contact \emph{cell}
decompositions. Although we do not delve into the details here, it
seems feasible to set up a more direct connection between the two
existence results just as in the topological case. The reader is
advised to turn to \cite{gs} or \cite{ma} for necessary background
on handle decompositions of manifolds and to \cite{e}, \cite{ge},
\cite{G} and \cite{ozst} for the related material on contact
topology.

\vspace{0.2in} \noindent{\bf {Acknowledgements.}} We are grateful
to John Etnyre, Hansj\"{o}rg Geiges, Andr\'{a}s Stipsicz, and
Emmanuel Giroux  for very useful comments on a draft of this
paper. We thank John Etnyre, in addition,  for sending us an
English translation of Giroux's paper {\em Convexit\'{e} en
topologie de contact,} by Daniel Mathews. We also thank Tolga
Etg\"{u} for helpful conversations. The author was partially
supported by the research grant 107T053 of the Scientific and
Technological Research Council of Turkey.

\section{Contact handles in dimension three}\label{handles}

We first review Giroux's contact handles in dimension three
\cite{gi}. The contact structure $\zeta_0 = \ker \a_0$, where
$\a_0 = dz-ydx + xdy$,  is the standard tight contact structure in
$\mathbb{R}^3$ and  the flow of the vector field
$$Z_0 = x \dfrac{\p}{\p x} + y\dfrac{\p}{\p y} + 2z \dfrac{\p}{\p
z}$$ preserves $\zeta_0$. Let $B^3 = \{ (x,y,z) \in \mathbb{R}^3
\;| \; x^2 + y^2 + z^2 \leq 1 \} $.  Then $\p B^3$ is a convex
surface since $Z_0$ is transverse to $\p B^3$. It is clear that
$Z_0$ lies in the contact planes $\zeta_0$ whenever $\a_0
(Z_0)=0$, i.e., when $z=0$. In other words, the disk $ B^3 \cap \{
z=0 \}$ is the characteristic surface in $B^3$ and $\p B^3 \cap \{
z=0 \}$ is the dividing curve on $\p B^3$.  \\

\begin{figure}[ht]
  \relabelbox \small {
  \centerline{\epsfbox{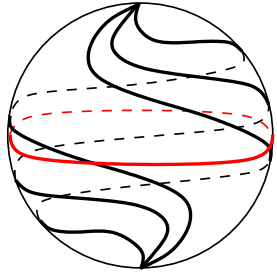}}}

  \endrelabelbox
        \caption{Characteristic foliation and the dividing curve on $\p B^3$ }
        \label{carfolhandle0}
\end{figure}

A model for a \emph{contact $0$-handle} is given as $(B^3,
\zeta_0)$, where $Z_0$ is used in gluing this handle.  Here the
orientation of the contact $0$-handle coincides with the usual
orientation of $\mathbb{R}^3$ (given by $dx \wedge dy \wedge dz$)
and its boundary has the induced orientation. The dividing curve
 divides the convex sphere $\p B^3$ into its
positive and negative regions: $R_+ = \p B^3 \cap \{ z > 0\}$  and
$R_- = \p B^3 \cap \{ z < 0\}$. The characteristic foliation on $\p B^3$ appears as
in Figure~\ref{carfolhandle0}, where the ``equator" is the dividing curve. \\

Giroux's criterion \cite{gi} implies that the dividing curve on
any tight $3$-ball with convex boundary is connected. Moreover
there is a unique tight contact structure on the $3$-ball with a
connected dividing set on its convex boundary up to isotopy fixing
the dividing set \cite{eli}. Hence we make the following
definition.

{\Def \label{stball} A standard contact $3$-ball is a tight
contact $3$-ball with convex boundary.} \\

As a matter of  fact, a contact $0$-handle is a model for a
standard contact $3$-ball and when we want to glue such a handle,
we use the vector field $Z_0$ in the model to obtain a ``contact"
collar neighborhood. A model for a  \emph{contact $3$-handle}, on
the other hand, is also defined as $(B^3, \zeta_0) $, where we
give opposite orientation to its boundary and use $-Z_0$ to glue
this
handle. \\

Let  $\zeta_1$ denote the contact structure in $\mathbb{R}^3$
given by the kernel of the $1$-form  $$\a_1 = dz + ydx + 2xdy,$$
and consider the vector field $$Z_1 = 2x \dfrac{\p}{\p x} -
y\dfrac{\p}{\p y} + z \dfrac{\p}{\p z}$$ whose flow preserves
$\zeta_1$. Observe that $\zeta_1$ is isotopic to the standard
tight contact structure $\zeta_0$ in $\mathbb{R}^3$.   Moreover,
for any $\e
> 0$,   $Z_1$ is transverse to the surfaces  $$ \{ (x,y,z) \in
\mathbb{R}^3 \; | \; x^2 + z^2 = \e^2 \} \; \mbox{and}  \; \{
(x,y,z) \in \mathbb{R}^3 \; | \; y^2 = 1 \}. $$ Note that the
intersection of these convex surfaces is not Legendrian. Let $$H_1
= \{ (x,y,z) \in \mathbb{R}^3 \; | \; x^2 + z^2 \leq \e^2, \; y^2
\leq 1 \} \; \mbox{and} \; F_1 = H_1 \cap \{ y=\pm 1 \}.$$ \\

\begin{figure}[ht]
  \relabelbox \small {
  \centerline{\epsfbox{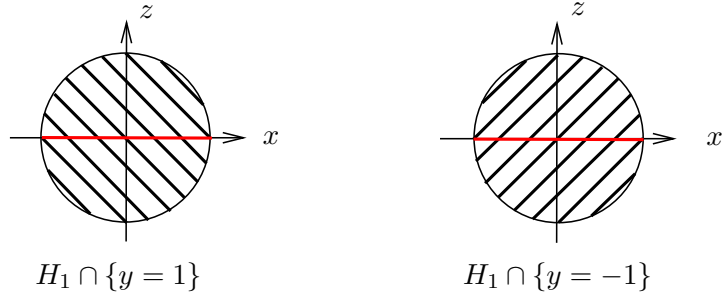}}}

\relabel{1}{{$x$}} \relabel{2}{{$z$}} \relabel{3}{{$x$}}
\relabel{4}{{$z$}} \relabel{a}{{$H_1 \cap \{ y = 1\} $}}
\relabel{b}{{ $H_1 \cap \{ y = -1\} $}}
  \endrelabelbox
        \caption{Characteristic foliation and the dividing set on $F_1 = H_1 \cap \{ y = \pm 1\} $}
        \label{carfolhandle1}
\end{figure}

A model for a \emph{contact $1$-handle} is given as $(H_1,
\zeta_1)$, where $Z_1$ is used in gluing this handle. Here the
contact $1$-handle acquires the usual orientation of
$\mathbb{R}^3$ and $\zeta_1$ orients $F_1$ as the outward pointing
normal vector field.  The characteristic surface in $H_1$ is given
by $H_1 \cap \{ z=0 \}$. The dividing curve $ \p H_1 \cap \{ z=0
\}$ divides $\p H_1$ into its positive and negative regions: $R_+
= \p H_1 \cap \{ z > 0\}$  and $R_- = \p H_1 \cap \{ z < 0\}$. The
characteristic foliation on $F_1$ is linear with slope $\mp 1$ on
 $H_1 \cap \{ y = \pm 1\} $ (viewed in a copy of the $xz$-plane) as depicted
 in Figure~\ref{carfolhandle1}. \\

A model for a \emph{contact $2$-handle} is defined as $(H_2,
\zeta_1)$, where $H_2 = \{ (x,y,z) \in \mathbb{R}^3 \; | \; x^2 +
z^2 \leq 1, \; y^2 \leq \e^2\} $. Note that the intersection of
the convex surfaces ($\e > 0$), $$ \{ (x,y,z) \in \mathbb{R}^3 \;
| \; x^2 + z^2 = 1\} \; \mbox{and} \; \{ (x,y,z) \in \mathbb{R}^3
\; |  \; y^2 = \e^2 \} $$ is not Legendrian. The characteristic
surface in $H_2$ is given by $H_2 \cap \{ z=0 \}$ and the dividing
curve on the boundary of the contact $2$-handle is given by $ \p
H_2 \cap \{ z=0 \}$. Let $F_2 = H_2 \cap \{ x^2+ z^2 = 1\} $. Here
the contact $2$-handle is oriented by the usual orientation of
$\mathbb{R}^3$; $-Z_1$ orients $F_2$ as the outward normal vector
field and we use $-Z_1$ when we glue such a handle along $F_2$. \\

\begin{figure}[ht]
  \relabelbox \small {
  \centerline{\epsfbox{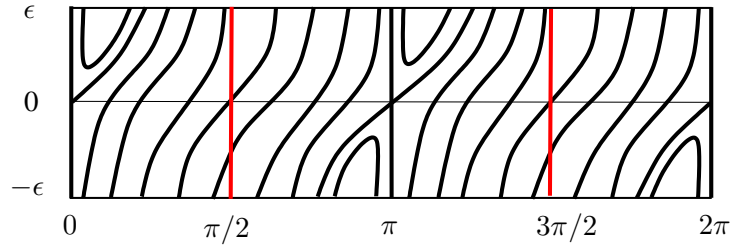}}}

\relabel{1}{{$0$}} \relabel{2}{{${\pi}/{2}$}} \relabel{3}{{$\pi$}}
\relabel{4}{{${3\pi}/{2}$}} \relabel{5}{{$2\pi$}}
\relabel{6}{{$\e$}} \relabel{7}{{0}} \relabel{8}{{$-\e$}}

  \endrelabelbox
        \caption{Characteristic foliation and the dividing set on $F_2$}
        \label{carfolhandle2}
\end{figure}

If we parametrize $F_2$ by $(\theta, y) \to (x= \sin \theta, y, z=
\cos \theta ) $ for  $(\theta,y) \in [0,2\pi] \times [-\e, \e]$,
then the equation for determining the characteristic foliation on
$F_2$ becomes $$(y\cos \theta - \sin \theta) d\theta + 2 \sin
\theta dy =0,$$ where the orientation of $F_2$ is given by
$d\theta \wedge dy$. Therefore the characteristic foliation is the
singular foliation which is given as the integral curves of the
equation $$\dfrac{dy}{d\theta}= \dfrac{1}{2} (1- y \cot \theta).$$
It follows that the characteristic foliation on $F_2$ appears as
in Figure~\ref{carfolhandle2}. Note that there are two hyperbolic
singular points corresponding to $(\theta, y) \in \{ (0, 0),(\pi,
0) \} $ and the dividing set on $F_2$ consists of the lines
$\theta=\pi /2$ and $\theta = 3 \pi /
2$. \\

Roughly speaking, a $3$-dimensional contact $k$-handle is a
topological $k$-handle which carries a tight contact structure
whose diving set on the boundary is depicted in
Figure~\ref{handles}. Moreover the characteristic foliations on
the gluing regions of these handles are shown in
Figures~\ref{carfolhandle0}, \ref{carfolhandle1} and
\ref{carfolhandle2}. \\

\begin{figure}[ht]
  \relabelbox \small {
  \centerline{\epsfbox{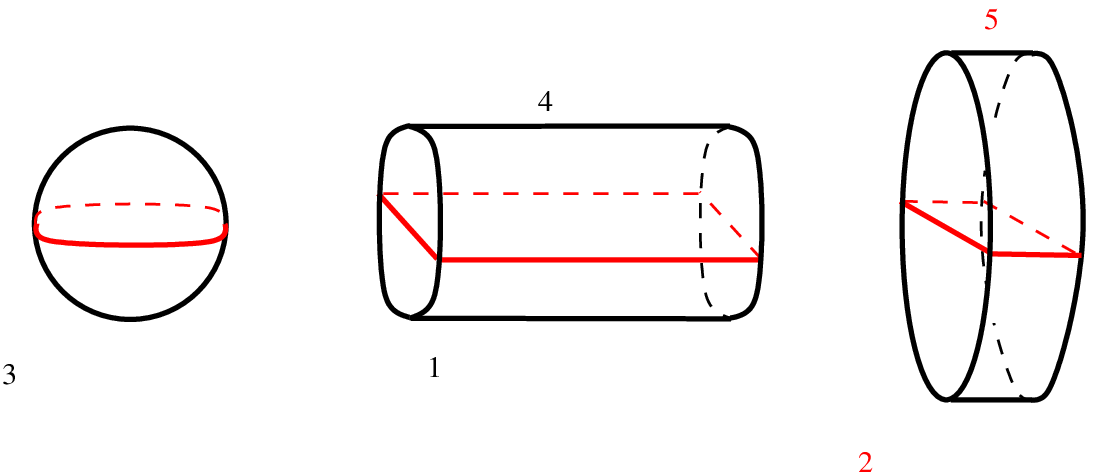}}}

\relabel{1}{{$1$-handle $D^1 \times D^2$}} \relabel{2}{{$2$-handle
$D^2 \times D^1$}} \relabel{3}{{$0$ and $3$-handle $D^3$}}
\relabel{4}{{$D^1$}} \relabel{5}{{$D^1$}}
  \endrelabelbox
        \caption{$3$-dimensional contact handles}
        \label{handles}
\end{figure}

Recall \cite{hon} that if two convex surfaces inside an ambient
contact $3$-manifold admit a Legendrian curve as their common
boundary, then the diving curves on these convex surfaces will
intersect that Legendrian curve in an ``alternating" fashion. In
the description of the contact $k$-handle, for $k=1,2$, however,
the diving curves on the convex surfaces which make up the
boundary of the contact $k$-handle do not meet the intersection of
these convex surfaces at an alternating fashion (see
Figure~\ref{handles}). This is not a contradiction since the
intersection of those convex surfaces is not Legendrian.

Next we would like to discuss contact handle attachments
\cite{gi}. By attaching contact $0$-handles we will  just mean
taking a disjoint union of some contact $0$-handles.  In order to
attach a contact $3$-handle to a contact $3$-manifold $(M, \xi)$
with convex boundary, we require that $\p M$ has at least one
component which is a $2$-sphere with a connected dividing set.
Then a contact $3$-handle attachment is just filling in this
$2$-sphere by a standard contact $3$-ball. The key point is that
the image of the characteristic foliation on the boundary $F_3 =
\p H_3$ of the $3$-handle under the attaching map is adapted to
the dividing curve $\G_{\p M}$ and therefore Giroux's Theorem
(\cite{gi}, Proposition II.3.6) allows us to glue the
corresponding contact structures.

Suppose that $(M,\xi)$ is a contact $3$-manifold with convex
boundary, where $\G_{\p M}$ denotes the dividing set on $\p M$. In
order to attach a contact $1$-handle to $M$ along two points $p$
and $q$ on $\G_{\p M}$ we identify the attaching region $F_1 \cong
D^0 \times D^2$ of the $1$-handle $H_1 \cong D^1 \times D^2$ with
regular neighborhoods of these points in $\p M$. The difference
from just a topological $1$-handle attachment is that we require
the dividing set on the attaching region  of the contact
$1$-handle to coincide with $\G_{\p M}$ on $\p M$ so that we can
glue the contact structures on $M $ and the contact $1$-handle
again by Giroux's Theorem (\cite{gi}, Proposition II.3.6). The
idea here is that once we initially identify the dividing curves
then we can match the characteristic foliations on the
\emph{convex} pieces that we glue by appropriate isotopies in the
collar neighborhoods given by the contact vector fields. Also we
need to make sure that the positive and the negative regions on
the corresponding convex boundaries match up so that the new
convex boundary after the handle attachment has well-defined $\pm$
regions divided by the new dividing set.

\begin{figure}[ht]
  \relabelbox \small {
  \centerline{\epsfbox{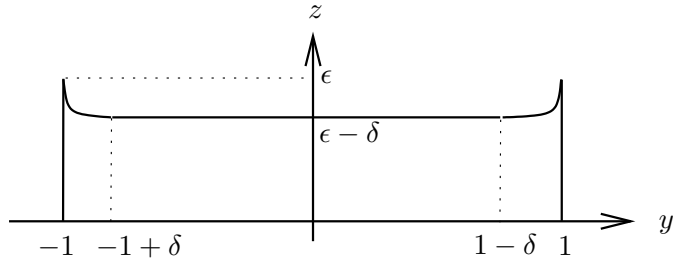}}}

\relabel{4}{{$1$}} \relabel{y}{{$y$}} \relabel{z}{{$z$}}

\relabel{1}{{$-1$}} \relabel{3}{{$1-\d$}} \relabel{2}{{$-1+\d$}}

\relabel{5}{{$\e$}} \relabel{6}{{$\e-\d$}}

  \endrelabelbox
        \caption{Modification of the contact $1$-handle}
        \label{rounding}
\end{figure}

Note that a contact $1$-handle is a manifold with \emph{corners}.
To get a smooth contact manifold as a result of a contact
$1$-handle attachment we propose the following modification
(similar to Honda's edge rounding technique \cite{hon}) to the
handle: Let $\d < \e$ be a sufficiently small positive real number
and let $f : [0,1] \to \mathbb{R}$ be a function defined as
follows:

\begin{itemize}
  \item $f(y) = \e-\d $ for $y \in [-1+\d, 1-\d]$,
  \item f is smooth on $(-1,1)$,
  \item  f is concave up on both $(-1, -1 + \d)$ and $(1-\d, 1)$,

  \item $\lim_{y \to \pm 1} f'(y) = \pm \infty$, and

\item $f (\pm 1)= \e$.
  \end{itemize}

\begin{figure}[ht]
  \relabelbox \small {
  \centerline{\epsfbox{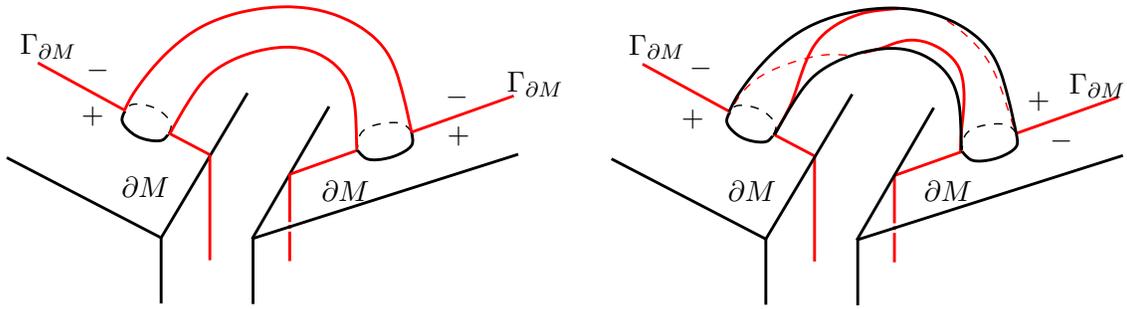}}}

\relabel{1}{{$\G_{\p M}$}} \relabel{4}{{$\p M$}}

\relabel{2}{{$\G_{\p M}$}} \relabel{3}{{$\p M$}}

\relabel{5}{{$+$}} \relabel{6}{{$-$}}

\relabel{7}{{$+$}} \relabel{8}{{$-$}}

\relabel{9}{{$\G_{\p M}$}} \relabel{12}{{$\p M$}}

\relabel{10}{{$\G_{\p M}$}} \relabel{11}{{$\p M$}}

\relabel{13}{{$+$}} \relabel{14}{{$-$}}

\relabel{15}{{$+$}} \relabel{16}{{$-$}}

\endrelabelbox
\caption{Attaching contact $1$-handles} \label{1-handle}
\end{figure}

Such a function is depicted in Figure~\ref{rounding}. Now consider
the region in the upper half $yz$-plane under the graph of the
function $f$ defined over $-1 \leq y \leq 1$. By revolving this
region around the $y$-axis, topologically we get a $1$-handle
(which looks like a vase).  One can verify that the contact vector
field $Z_1$ is still transverse to the side surface as well as the
top and the bottom disks. When we glue this (modified) contact
$1$-handle to a contact manifold with convex boundary we get a
smooth manifold carrying a contact structure which makes the
resulting boundary convex. In Figure~\ref{1-handle}, we
illustrated two possible contact $1$-handle attachments (taking
into account the compatibility of the $\pm$ regions),  where
corners should be smoothed as explained above.

\begin{figure}[ht]
  \relabelbox \small {
  \centerline{\epsfbox{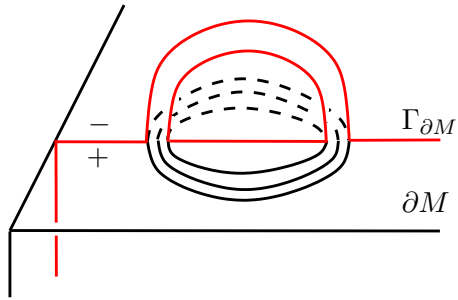}}}

\relabel{4}{{$\p M$}}

\relabel{1}{{$\G_{\p M}$}}

\relabel{3}{{$+$}} \relabel{2}{{$-$}}

  \endrelabelbox
        \caption{Attaching a contact $2$-handle}
        \label{2-handle}
\end{figure}

Next we explain how to attach a contact $2$-handle on top of a
contact $3$-manifold $(M, \xi)$ with convex boundary. As expected,
attachment of a contact $2$-handle requires more work compared to
the other contact handles. The attaching curve is the image of the
core circle of the annulus $F_2 \cong \p D^2 \times D^1$ under the
attaching map $F_2 \to \p M$ of the $2$-handle. It is well-known
that in order to attach a topological $2$-handle one only has to
specify the attaching curve on $\p M$. To attach a contact
$2$-handle, however, we require the attaching curve to intersect
$\G_{\p M}$ transversely at two distinct points. This will allow
one to glue the contact structures on $M$ and the $2$-handle as
explained in great details by Giroux in (\cite{gi}, Lemma
III.3.2). The idea here is that one can construct a singular
foliation adapted to $\G_{\p M}$ which conjugates to the
characteristic foliation on $F_2$ (see Figure~\ref{carfolhandle2})
in an annulus neighborhood of the attaching curve on the convex
surface $\p M$. In addition, just as in attaching a contact
$1$-handle, we need to pay attention so that the $\pm$ regions in
the corresponding boundaries match up appropriately. Moreover, one
can smooth the corners of the contact $2$-handle by a modification
which preserves the convexity of the boundary---similar to the
modification we explained for contact $1$-handles.

\section{Contact handle decompositions}\label{decomposition}

{\Thm [Giroux] Every compact contact $3$-manifold (closed or with
convex boundary) admits a contact handle decomposition.}

\begin{proof} Suppose that $(M, \xi)$ is a \emph{connected} contact
$3$-manifold with convex boundary.  It follows that $(M, \G_{\p
M}, \xi)$ admits a compatible partial open book decomposition
\cite{hkm} and, in particular,  $(M, \xi)$ can be decomposed into
two tight contact handlebodies $(H, \xi\vert_{ H})$ and $(N,
\xi\vert_{ N})$ where $H$ is connected by our assumption that $M$
is connected (see \cite{eo} for notation). Now we claim that $(H,
\xi\vert_{ H})$ has a contact handle decomposition with a unique
contact $0$-handle and some contact $1$-handles. This is because
$(H, \xi\vert_{ H})$ is \emph{product disk decomposable}
\cite{hkm}, i.e., there exist some pairwise disjoint compressing
disks in $H$ each of whose boundary intersects $\G_{\p H}$
transversely in two points,  so that when we cut $H$ along these
disks we get a standard contact $3$-ball. Clearly the resulting
standard contact $3$-ball can be considered as a contact
$0$-handle. On the other hand, the thickening of a compressing
disk satisfies our definition of a contact $1$-handle which is
attached to the contact $0$-handle. This proves our claim about
the tight contact handlebody $(H, \xi\vert_{ H})$. Moreover each
component of the handlebody $N$ is also product disk decomposable.
By turning the handles upside down we conclude that $(M,
\xi\vert_M)$ is obtained from $(H, \xi\vert_H)$ by attaching some
contact $2$ and $3$-handles. Thus we proved that $(M, \xi)$ admits
a contact handle decomposition.

Suppose that $(Y, \xi)$ is a \emph{connected} and \emph{closed}
contact $3$-manifold. Let $p$ be an arbitrary point in $Y$. Then,
by Darboux's theorem, there is a neighborhood of $p$ in $Y$ which
is just a standard contact $3$-ball. Now the closure of the
complement of this ball in $Y$ is a contact $3$-manifold $(M,
\xi\vert_M)$ whose boundary is a convex $2$-sphere with a
connected dividing set $\G_{\p M}$. We proved above that $(M,
\xi)$ admits a contact handle decomposition. Furthermore we can
obtain $(Y, \xi)$ from $(M, \xi\vert_M)$ by gluing back the
standard contact $3$-ball that we deleted at the beginning, which
is indeed equivalent to a contact $3$-handle attachment. Hence we
proved that $(Y, \xi)$ has a contact handle decomposition with a
unique contact $0$-handle and some contact $1,2$ and $3$-handles.
If $(Y, \xi)$ is not connected then we can apply the above
argument to each of its components to obtain a contact handle
decomposition.

\end{proof}

\section{Bypass attachment}\label{bypass}

Recall that a bypass \cite{hon} for a convex surface $\S$ in a
contact $3$-manifold is an oriented embedded half-disk $D$ with
Legendrian boundary, satisfying the following:

 \begin{itemize}

 \item  $\p D $ is the union of two arcs $\g_1$ and $\g_2$ which
intersect at their endpoints,

\item $D$ intersects $\S$ transversely along $\g_2$,

\item $D$ (or $D$ with the opposite orientation) has the following
tangencies along $\p D$:

(1) positive elliptic tangencies at the endpoints of $\g_2$,

(2) one negative elliptic tangency on the interior of $\g_2$, and

(3) only positive tangencies along $\g_1$, alternating between
elliptic and hyperbolic,

\item $\g_2$ intersects the dividing set $\G$ exactly at three
points, and these three points are the elliptic points of $\g_2$.

 \end{itemize}

\begin{figure}[ht]
  \relabelbox \small {
  \centerline{\epsfbox{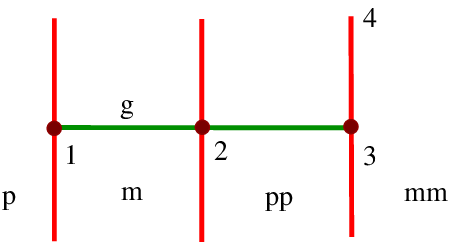}}}

 \relabel{1}{{$p_1$}}  \relabel{2}{{$p_2$}}  \relabel{3}{{$p_3$}}
\relabel{g}{{$\gamma_2$}}
 \relabel{4}{{$\Gamma$}} \relabel{p}{{$-$}} \relabel{m}{{$+$}}
 \relabel{pp}{{$-$}} \relabel{mm}{{$+$}}

  \endrelabelbox
        \caption{The attaching arc of a bypass intersecting the dividing
        set $\Gamma$ at $\{p_1, p_2, p_3\}$}
        \label{bypass-0}
\end{figure}

In this section we show that a bypass attachment consists of a
pair of \emph{contact} 1 and 2-handles---which cancel each other
out only topologically. Here by a bypass attachment we mean
attaching a thickened neighborhood of the bypass disk $D$. The
attaching arc $\g_2$ of an \emph{exterior} bypass is a Legendrian
arc on the convex boundary of a contact $3$-manifold, where $\g_2$
intersects $\G$ transversely at $p_2$, and $p_1, p_3 \in \G$ are
the endpoints of $\g_2$, as we depict in Figure~\ref{bypass-0}.

\begin{figure}[ht]
  \relabelbox \small {
  \centerline{\epsfbox{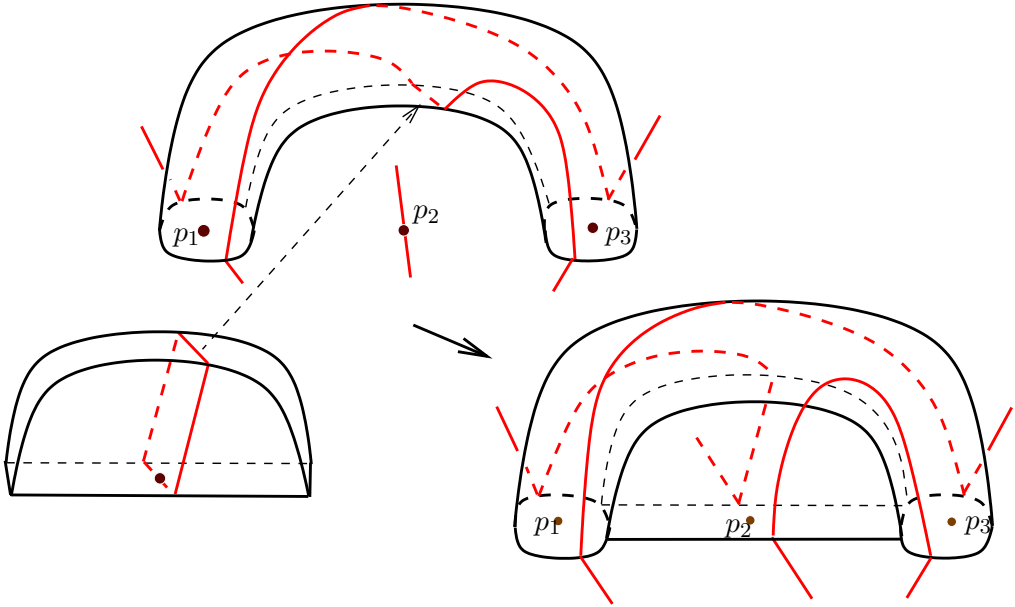}}}

\relabel{1}{{$p_1$}}    \relabel{3}{{$p_3$}} \relabel{2}{{$p_2$}}

\relabel{11}{{$p_1$}}    \relabel{33}{{$p_3$}}
 \relabel{22}{{$p_2$}}

  \endrelabelbox
        \caption{Anatomy of a bypass attachment}
        \label{bypass-1}
\end{figure}

In order to attach a bypass along the arc $\g_2$ indicated in
Figure~\ref{bypass-0}, we first attach a contact 1-handle  whose
feet are identified with the neighborhoods of $p_1$ and $p_3$,
respectively. Here we pay attention to the compatibility of the
$\pm$ regions in the surfaces that we glue together. To be more
precise, we describe the gluing map $\phi$ which identifies the
gluing region $F_1$ of the contact $1$-handle with two disjoint
disks around $p_1$ and $p_3$ as follows: $\phi$ takes $(0,-1,0)$
to $p_1$, the dividing arc $\{ -1 \leq x \leq 1, y=-1, z=0 \}$ to
an arc around $p_1$ in $\G$ and the arc $\{ x=0, y=-1,  -1 \leq z
\leq 0\}$ to an arc on $\g_2$. Similarly, $\phi$ takes $(0,1,0)$
to $p_3$, the dividing arc $\{ -1 \leq x \leq 1, y=1, z=0 \}$ to
an arc around $p_3$ in $\G$ and the arc $\{ x=0, y=1, 0 \leq z
\leq 1 \}$ to an arc on $\g_2$ (see Figure~\ref{bypass-1}). Now we
claim that we can attach a topologically cancelling contact
$2$-handle so that the union of the contact $1$ and $2$-handles
that we attach has the same effect as attaching  a bypass along
$\g_2$. Hence this procedure gives the contact anatomy of a bypass
attachment, which is depicted (locally) in Figure~\ref{bypass-1}.
In the following, we explain how to glue the contact $2$-handle so
that the union of the contact $1$ and $2$-handles can be viewed as
a neighborhood of a bypass disk $D =D_1 \cup D_2$, where $D_i$
is a disk in the contact $i$-handle, for $i=1,2$.\\

\begin{figure}[ht]
  \relabelbox \small {
  \centerline{\epsfbox{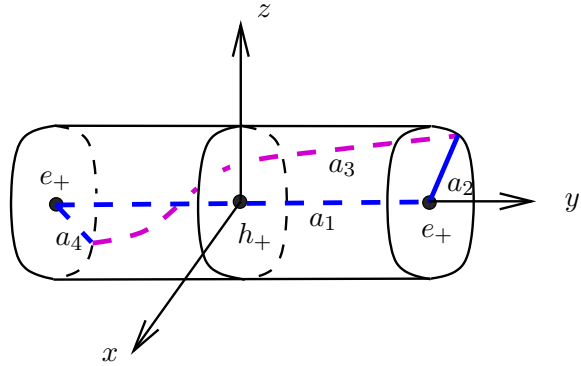}}}

\relabel{1}{{$a_1$}}  \relabel{2}{{$a_2$}}  \relabel{3}{{$a_4$}}
\relabel{4}{{$a_3$}} \relabel{x}{{$x$}} \relabel{y}{{$y$}}
\relabel{z}{{$z$}} \relabel{e}{{$e_+$}}  \relabel{ee}{{$e_+$}}
 \relabel{h}{{$h_+$}}
  \endrelabelbox
        \caption{The disk $ D_1$ has boundary $a_1 \cup a_2 \cup a_3 \cup a_4$. }
                \label{convex1}
\end{figure}

Construction of $D_1$: The idea here is to perturb the
(rectangular) disk $\{z=0\} \cap \{ x \leq 0 \} $ in the contact
$1$-handle $(H_1, \zeta_1) $ so that the boundary of that disk is
a ``Legendrian" curve on which there are one positive hyperbolic
and two positive elliptic singular points. To be more precise, let
$a_1$ denote the Legendrian arc $\{x=z=0\}$ in $H_1$; $a_2$ denote
the Legendrian arc $\{y=1\} \cap \{x=-z\} \cap \{ z \geq 0\}$;
$a_3$ denote a Legendrian arc connecting the points
 $(-\dfrac{\e}{\sqrt{2}}, 1, \dfrac{\e}{\sqrt{2}})$ and
$(-\dfrac{\e}{\sqrt{2}}, -1, -\dfrac{\e}{\sqrt{2}})$ on $\p H_1$
(see Figure~\ref{convex1}) and $a_4$ denote the Legendrian arc
$\{y=-1\} \cap \{x=z\} \cap \{z \leq 0\}$. Then $a_1 \cup a_2 \cup
a_3 \cup a_4$ bounds a surface $D_1$ in $(H_1, \zeta_1) $, where
$(0,0,0)$ is a hyperbolic singular point and $(0, \pm 1, 0)$ are
elliptic singular points on $\p D_1$. Moreover we orient $D_1$
such that all the
singularities on $\p D_1$ are positive. \\

\begin{figure}[ht]
  \relabelbox \small {
  \centerline{\epsfbox{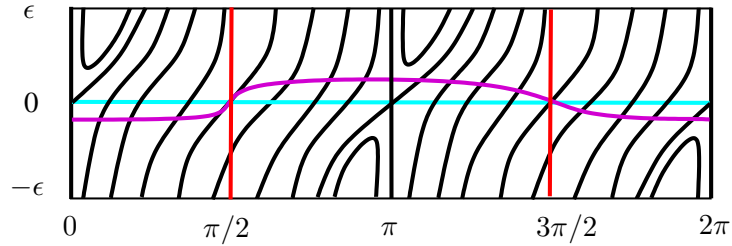}}}

\relabel{1}{{$0$}} \relabel{2}{{${\pi}/{2}$}} \relabel{3}{{$\pi$}}
\relabel{4}{{${3\pi}/{2}$}} \relabel{5}{{$2\pi$}}
\relabel{6}{{$\e$}} \relabel{7}{{0}} \relabel{8}{{$-\e$}}

  \endrelabelbox
        \caption{Perturbation of the curve $\{y=0 \} \cap F_2$}
        \label{diskinhandle2}
\end{figure}

Construction of $D_2$: The idea here is to perturb the disk
$\{y=0\} \cap H_2$ in the contact $2$-handle $(H_2, \zeta_1) $
into a disk whose boundary is a Legendrian circle on which there
is a unique elliptic singularity. To achieve this we first perturb
the curve $\{y=0\}$ on $F_2$ as follows: Fix the points $(\theta,
y) \in \{ (\pi/2, 0), (3 \pi/2, 0)\}$ and push the arc $\{\pi/2
\leq \theta \leq 3 \pi/2 , y=0 \}$ slightly in the upward
direction and the arc $\{ 0 \leq \theta \leq \pi/2, y=0 \} \cup
\{3 \pi/2 \leq \theta \leq 2\pi, y=0\}$ slightly in the downward
direction as shown in Figure~\ref{diskinhandle2}; Legendrian
realize the perturbed curve and then consider the spanning disk
$D_2$. With a little bit of care, we can make sure that $\p D_2$
has a unique elliptic singular point at $\theta= \pi/2$. More
precisely, to have an elliptic singularity at $\theta= \pi/2$, the
slope of the perturbed curve should agree with the slope of the
characteristic foliation at that point
on $F_2$, which certainly can be arranged. \\

\begin{figure}[ht]
  \relabelbox \small {
  \centerline{\epsfbox{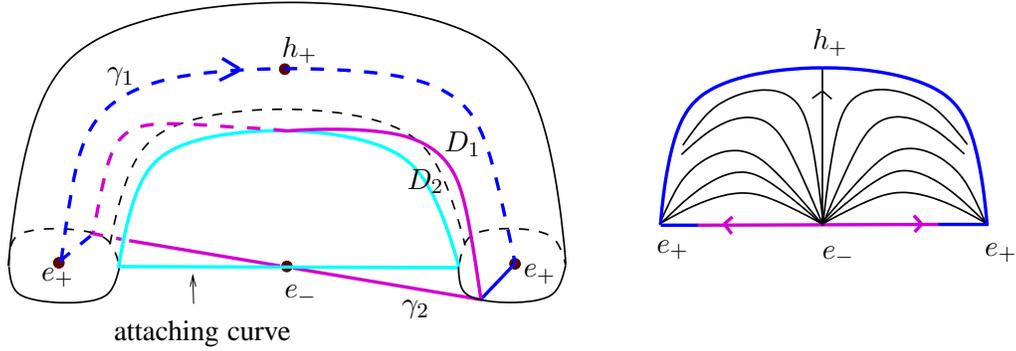}}}

\relabel{1}{{$\g_1$}}  \relabel{2}{{$\g_2$}}  \relabel{3}{{$D_2$}}
\relabel{4}{{$D_1$}} \relabel{5}{{$e_+$}} \relabel{6}{{$e_+$}}
\relabel{7}{{$e_-$}} \relabel{8}{{$h_+$}} \relabel{a}{{$e_+$}}
\relabel{c}{{$e_+$}} \relabel{b}{{$e_-$}} \relabel{d}{{$h_+$}}
\relabel{cu}{{attaching curve}}

  \endrelabelbox
        \caption{Left: Bypass disk $D= D_1 \cup D_2$ inside a bypass attachment; Right:
        The characteristic foliation on $D$}
                \label{bydisk}
\end{figure}

In order to exhibit the bypass disk $D$, we glue the disks $D_1$
and $D_2$ along some parts of their boundaries as follows. Let us
express $\p D_2$ as a union of two arcs $b_1$ and $b_2$ where $b_1
= \p D_2 \cap \{0 \leq \theta \leq \pi \}$ and $b_2 = \p D_2 \cap
\{\pi \leq \theta \leq 2\pi \}$ on $F_2$. Then $D$ is obtained by
gluing $D_1$ and $D_2$ where we simply identify $a_3$ and $b_2$.
This can be achieved if the attaching diffeomorphism takes the
core $\{ y =0 \}$ of the attaching region $F_2$ of the 2-handle
$H_2$ to the attaching curve that is indicated in
Figure~\ref{bydisk}. Note that the boundary of the disk $D$
consists of the Legendrian arcs $\g_1 = a_1$ and $\g_2 = a_2 \cup
b_1 \cup a_4$ and hence we can isotope $D$ to be convex.  If we
orient $D$ keeping the orientation of $D_1$, then the sign of the
unique elliptic point on $\p D_2$ becomes negative. The
characteristic foliation on the convex disk $D$ appears as in
Figure~\ref{bydisk}, since the negative elliptic singular point is
a source whereas the positive elliptic singular points are sinks
and there is a unique
hyperbolic singular point on $\p D$.  \\

\section{An infinite family of overtwisted contact $3$-spheres}\label{examples}

\begin{figure}[ht]
  \relabelbox \small {
  \centerline{\epsfbox{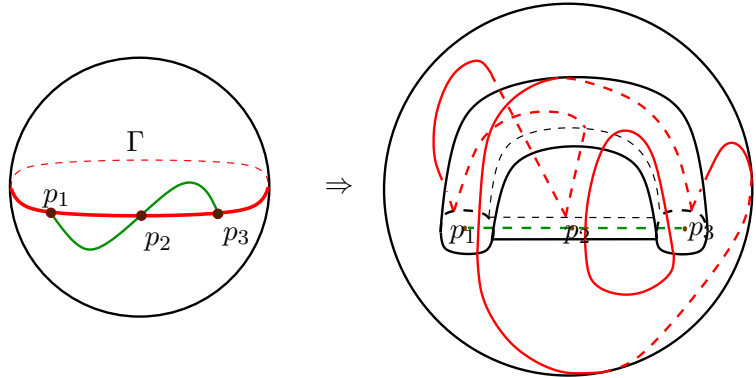}}}

 \relabel{1}{{$p_1$}}  \relabel{2}{{$p_2$}}  \relabel{3}{{$p_3$}}
\relabel{11}{{$p_1$}}    \relabel{33}{{$p_3$}}
\relabel{r}{{$\Rightarrow$}}  \relabel{22}{{$p_2$}}

 \relabel{4}{{$\Gamma$}}

  \endrelabelbox
        \caption{The result of a bypass attachment to a contact $0$-handle}
        \label{bypass-arc-2}
\end{figure}

\noindent \textbf{An overtwisted contact $3$-sphere:} In the
following we describe a contact handle decomposition of an
overtwisted contact structure $\xi_{0}$ in $S^3$. We start with
attaching a bypass to a contact $0$-handle along the Legendrian
arc depicted in Figure~\ref{bypass-arc-2} on the convex sphere $\p
B^3$, where the \emph{southern} hemisphere is the $+$ region.

\begin{figure}[ht]
  \relabelbox \small {
  \centerline{\epsfbox{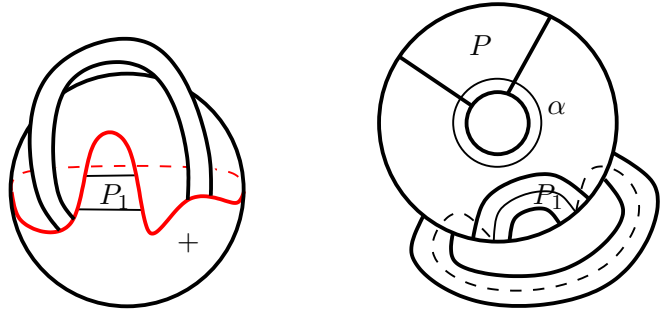}}}
 \relabel{1}{{$P_1$}}  \relabel{2}{{$+$}}

 \relabel{t}{{$P$}}  \relabel{P}{{$P_1$}}  \relabel{a}{{$\a$}}

  \endrelabelbox
        \caption{Left: The new page $S'$ is $S \;\cup$ the
        attached 1-handle and $P'=P \cup P_1$; Right: The $1$-handle $P$
        and a right-handed Dehn twist around $\a$ can be
        viewed as a stabilization of the rest.}
        \label{partial-stabilization}
\end{figure}

The Legendrian arc has its endpoints  at $p_1, p_3 \in \G$ and
intersects $\G$ transversely at $p_2$. The diving set on the
convex boundary of the resulting $3$-ball $B^3_{ot}$ after the
bypass attachment consists of three connected components (see
Figure~\ref{bypass-arc-2}) and it follows, by Giroux's criterion
\cite{gi}, that $B^3_{ot}$ is overtwisted. Moreover we claim that
$B^3_{ot}$ is the standard neighborhood of an overtwisted disk. To
prove our claim we first describe a partial open book of
$B^3_{ot}$. The partial open book of a contact $0$-handle is
described in \cite{eo1}. The page $S$ is an annulus, $P$ is a
neighborhood of a trivial arc connecting the distinct components
of the boundary of this annulus, and the monodromy is a
right-handed Dehn twist along the core of the annulus. In
\cite{hkm}, Honda, Kazez and Mati\'{c} describe how to obtain a
partial open book of the resulting contact $3$-manifold after a
bypass attachment. According to their recipe, a 1-handle is
attached to the page $S$ to obtain the new page $S'$ as depicted
in Figure~\ref{partial-stabilization}. Moreover $P' = P \cup P_1$,
and the embedding of the new piece $P_1$ into $S'$ is described
explicitly in Figure~\ref{partial-stabilization}: The solid arc in
$P_1$ is mapped to the dashed arc going once over the new
$1$-handle. It follows that when we attach a bypass to a contact
$0$-handle along the arc given in Figure~\ref{bypass-arc-2}, the
resulting partial open book (see
Figure~\ref{partial-stabilization})  is nothing but a positive
stabilization of the partial open book of a standard neighborhood
of an overtwisted disk (\cite{hkm}, see also \cite{eo1}).

\begin{figure}[ht]
  \relabelbox \small {
  \centerline{\epsfbox{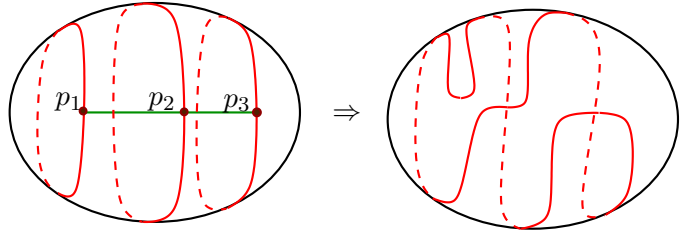}}}

 \relabel{1}{{$p_1$}}  \relabel{2}{{$p_2$}}  \relabel{3}{{$p_3$}}

\relabel{r}{{$\Rightarrow$}}

  \endrelabelbox
        \caption{Left: The attaching arc of a second bypass; Right:
        The dividing set after the second bypass attachment}
        \label{bypass-arc-3}
\end{figure}

Next we attach another bypass to $B^3_{ot}$ along the given arc on
$\p B^3_{ot}$ as depicted in Figure~\ref{bypass-arc-3}. The diving
set on the convex boundary of the resulting $3$-ball is connected
as shown in Figure~\ref{bypass-arc-3} and therefore we can cap off
the convex boundary by a contact $3$-handle. The resulting contact
$3$-sphere $(S^3,\xi_0)$, which consists of a contact $0$-handle,
two contact $1$-handles, two contact $2$-handles and a contact
$3$-handle,  is indeed  overtwisted. In fact, we will show that
$d_3 (\xi_{0}) = 1/2 $, which determines the homotopy (and hence
the isotopy) class of the overtwisted contact structure  $\xi_{0}$
in $S^3$.

To prove our claim we observe that

{\Rem \label{upsidedown} We can turn contact handles upside down
and a contact $k$-handle becomes a contact $(3-k)$-handle when
turned upside down. Moreover, a bypass turned upside down is
another bypass attached from the other side.} \\

Thus the second bypass and the last contact $3$-handle attached to
$B^3_{ot}$ can be viewed as a copy of $B^3_{ot}$ when the contact
handles are turned upside down. This is because the upside down
bypass is attached to the contact $0$-handle along an arc isotopic
to the one in Figure~\ref{bypass-arc-2}. Hence we conclude that
$(S^3, \xi_0)$ can be obtained by taking the double of the
standard neighborhood $B^3_{ot}$ of the overtwisted disk instead
of attaching the second bypass and the last contact $3$-handle.
Since we know a partial open book for  $B^3_{ot}$, we can actually
construct an open book for the double by ``gluing" the partial
open books along their boundaries as explained in \cite{hkm}. It
turns out \cite{eo1} that the open book for $(S^3, \xi_0)$ has a
twice punctured disk as its page and the monodromy is given by a
positive and a negative Dehn twists along the two punctures,
respectively. It is known (see, for example \cite{eto}) that the
$d_3$-invariant of the contact structure corresponding to such an
open book is equal to $1/2$. \\

\begin{figure}[ht]
  \relabelbox \small {
  \centerline{\epsfbox{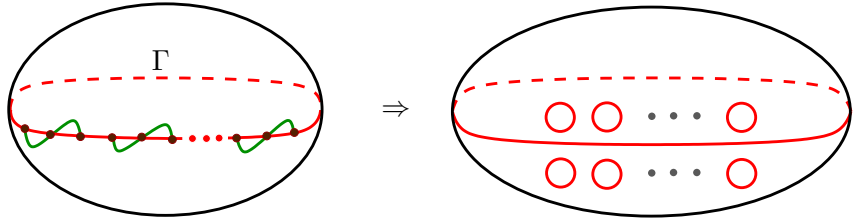}}}
 \relabel{4}{{$\Gamma$}} \relabel{r}{{$\Rightarrow$}}
 \endrelabelbox
        \caption{Left: The attaching arcs for $n$ bypasses;
        Right: The dividing set after attaching bypasses along the given arcs}
        \label{multi}
\end{figure}

\noindent \textbf{An infinite family of overtwisted contact
$3$-spheres:} We can generalize our discussion above to obtain
contact handle decompositions of infinitely many pairwise
non-isotopic overtwisted contact $3$-spheres. We first fix a
positive integer $n$, and choose a sequence of nearby points $p_1,
p_2, \ldots, p_{3n}$ on the dividing set on the boundary $\p B^3$
of the contact $0$-handle, where the \emph{southern}  hemisphere
is the $+$ region. For $k= 1,4,7, \ldots, 3n-2$, let $\gamma_k$ be
an arc isotopic to the one depicted in Figure~\ref{bypass-arc-2}
starting at $p_k$, passing through $p_{k+1}$, and ending at
$p_{k+2}$. Next we attach a bypass along each $\gamma_k$ to this
contact $0$-handle. The result of attaching these bypasses is
indeed an overtwisted $3$-ball where the dividing set on the
convex boundary has $2n+1$ connected components as shown in
Figure~\ref{multi}.

\begin{figure}[ht]
  \relabelbox \small {
  \centerline{\epsfbox{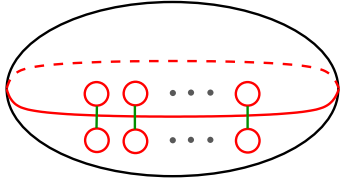}}}

 \endrelabelbox
        \caption{The attaching arcs for the second set of bypasses}
        \label{nby}
\end{figure}

The resulting partial open book can be constructed similar to the
$n=1$ case (that we already discussed), since a bypass attachment
is just a local modification. Then by taking the double of the
resulting overtwisted $3$-ball we get an overtwisted $3$-sphere
$(S^3, \xi_n)$. The page of the open book compatible with $(S^3,
\xi_n)$ is a disk with $2n$-punctures. Let  $t_{m}$ denote a
right-handed Dehn twist around $\a_m$, where  $\a_m$ is a curve
around the $m$th puncture. Then the monodromy of this open book is
given by $\prod_{i=1}^{n} t_{i} t^{-1}_{n+i} $. It follows that,
$d_3 (\xi_n)=(2n+1)/2 $, since $\xi_n$ can be obtained from
$\xi_{n-1}$ by a positive stabilization followed by a negative
stabilization, where a negative stabilization increases the
$d_3$-invariant by one while a positive stabilization does not
affect the contact structure (see, for example, \cite{ozst}).
Similar to the $n=1$ case, instead of doubling the overtwisted
$3$-ball to obtain $(S^3, \xi_n)$, we can attach $n$ more bypasses
to this $3$-ball along the $n$ arcs shown in Figure~\ref{nby} and
a contact $3$-handle to cap off the resulting boundary (see
Remark~\ref{upsidedown}). Hence, for each positive integer $n$, we
get an explicit contact handle decomposition of the overtwisted
$3$-sphere $(S^3, \xi_n)$ consisting of a contact $0$-handle, $2n$
contact $1$-handles, $2n$ contact $2$-handles and a contact
$3$-handle.

\section{Final Remarks}

It is well-known that one can slide handles in a given handle
decomposition of a smooth manifold. The natural question which
arises from the discussion in this paper is that whether there is
an analogue of handle sliding in contact topology. Similarly one
can ask whether there is a contact handle cancellation? It seems
to us that both questions have affirmative answers and we are
planning to investigate such issues in a future work.

In addition, it may be possible to compute the $EH$-class of a
contact $3$-manifold via its contact handle decomposition. In
order to achieve this goal one can first obtain a partial open
book decomposition of the contact $3$-manifold based on its handle
decomposition. The idea here is that the page $S$ of a partial
open book will acquire a $1$-handle once we attach a contact
$1$-handle to the contact $3$-manifold at hand. The attachment of
a contact $2$-handle (in fact just its attaching curve) will
simply determine $P$ and its embedding in $S$. The attachments of
contact $0$ and $3$-handles will manifest themselves merely as
suitable  stabilizations. Finally, to compute the $EH$-class of
the resulting contact $3$-manifold, we apply the techniques
recently developed by Honda, Kazez and Mati\'{c} \cite{hkm}.


\end{document}